\documentclass[12pt]{amsart}
\usepackage{amssymb,latexsym}
\usepackage{pdfsync}

\newdimen\AAdi%
\newbox\AAbo%
\def\AAk#1#2{\s_etbox\AAbo=\hbox{#2}\AAdi=\wd\AAbo\kern#1\AAdi{}}%
\def\AAr#1#2#3{\s_etbox\AAbo=\hbox{#2}\AAdi=\ht\AAbo\raise#1\AAdi\hbox{#3}}%
\font\tenmsb=msbm10 at 12pt
\font\sevenmsb=msbm7 at 8pt
\font\fivemsb=msbm5 at 6pt
\newfam\msbfam
\textfont\msbfam=\tenmsb
\scriptfont\msbfam=\sevenmsb
\scriptscriptfont\msbfam=\fivemsb
\def\Bbb#1{{\tenmsb\fam\msbfam#1}}

\newcommand{\beq}{\begin{equation}}
\newcommand{\eeq}{\end{equation}}
\newcommand{\beqr}{\begin{eqnarray}}
\newcommand{\eeqr}{\end{eqnarray}}
\newcommand{\ba}{\begin{array}}
\newcommand{\ea}{\end{array}}

\begin{document}

\newtheorem{thm}{Theorem}
\newtheorem{lem}{Lemma}
\newtheorem{cor}{Corollary}
\newtheorem{rem}{Remark}
\newtheorem{pro}{Proposition}
\newtheorem{defi}{Definition}
\newtheorem{eg}{Example}
\newtheorem*{claim}{Claim}
\newtheorem{conj}[thm]{Conjecture}
\newcommand{\noi}{\noindent}
\newcommand{\dis}{\displaystyle}
\newcommand{\mint}{-\!\!\!\!\!\!\int}
\numberwithin{equation}{section}

\def \bx{\hspace{2.5mm}\rule{2.5mm}{2.5mm}}
\def \vs{\vspace*{0.2cm}}
\def\hs{\hspace*{0.6cm}}
\def \ds{\displaystyle}
\def \p{\partial}
\def \O{\Omega}
\def \o{\omega}
\def \b{\beta}
\def \m{\mu}
\def \l{\lambda}
\def\L{\Lambda}
\def \ul{u_\lambda}
\def \D{\Delta}
\def \d{\delta}
\def \k{\kappa}
\def \s{\sigma}
\def \e{\varepsilon}
\def \a{\alpha}
\def \tf{\tilde{f}}
\def\cqfd{%
\mbox{ }%
\nolinebreak%
\hfill%
\rule{2mm} {2mm}%
\medbreak%
\par%
}
\def \pr {\noindent {\it Proof.} }
\def \rmk {\noindent {\it Remark} }
\def \esp {\hspace{4mm}}
\def \dsp {\hspace{2mm}}
\def \ssp {\hspace{1mm}}

\def\la{\langle}\def\ra{\rangle}

\def \u{u_+^{p^*}}
\def \ui{(u_+)^{p^*+1}}
\def \ul{(u^k)_+^{p^*}}
\def \energy{\int_{\R^n}\u }
\def \sk{\s_k}
\def \mo{\mu_k}
\def\cal{\mathcal}
\def \I{{\cal I}}
\def \J{{\cal J}}
\def \K{{\cal K}}
\def \OM{\overline{M}}

\def\n{\nabla}

\def\fk{{{\cal F}}_k}
\def\M1{{{\cal M}}_1}
\def\Fk{{\cal F}_k}
\def\Fl{{\cal F}_l}
\def\FF{\cal F}
\def\Gk{{\Gamma_k^+}}
\def\n{\nabla}
\def\uuu{{\n ^2 u+du\otimes du-\frac {|\n u|^2} 2 g_0+S_{g_0}}}
\def\uuug{{\n ^2 u+du\otimes du-\frac {|\n u|^2} 2 g+S_{g}}}
\def\sku{\sk\left(\uuu\right)}
\def\qed{\cqfd}
\def\vvv{{\frac{\n ^2 v} v -\frac {|\n v|^2} {2v^2} g_0+S_{g_0}}}
\def\vvs{{\frac{\n ^2 \tilde v} {\tilde v}
 -\frac {|\n \tilde v|^2} {2\tilde v^2} g_{S^n}+S_{g_{S^n}}}}
\def\skv{\sk\left(\vvv\right)}
\def\tr{\hbox{tr}}
\def\pO{\partial \Omega}
\def\dist{\hbox{dist}}
\def\RR{\Bbb R}\def\R{\Bbb R}
\def\C{\Bbb C}
\def\B{\Bbb B}
\def\N{\Bbb N}
\def\Q{\Bbb Q}
\def\Z{\Bbb Z}
\def\PP{\Bbb P}
\def\EE{\Bbb E}
\def\F{\Bbb F}
\def\G{\Bbb G}
\def\H{\Bbb H}
\def\SS{\Bbb S}\def\S{\Bbb S}

\def\div{\hbox{div}\,}

\def\lcf{{locally conformally flat} }

\def\circledwedge{\setbox0=\hbox{$\bigcirc$}\relax \mathbin {\hbox
to0pt{\raise.5pt\hbox to\wd0{\hfil $\wedge$\hfil}\hss}\box0 }}

\def\sss{\frac{\s_2}{\s_1}}

\date{December 18, 2024}
\title[ Rigidity theorems of Lagrangian and symplectic translators  ]{ Rigidity theorems of Lagrangian and symplectic translating solitons}

\author{}

\author[Qiu]{Hongbing Qiu}
\address{School of Mathematics and Statistics\\ Wuhan University\\Wuhan 430072,
China
 }
 \email{hbqiu@whu.edu.cn}

\begin{abstract}

By carrying out refined point-wise estimates for the mean curvature, we prove better rigidity theorems of Lagrangian and symplectic translating solitons.

\vskip12pt

\noindent{\it Keywords and phrases}:  Lagrangian, symplectic, translating soliton, rigidity, mean curvature

\noindent {\it MSC 2020}:  53C24, 53E10 

\end{abstract}
\maketitle
\section{Introduction}

Let $J, \omega$ be the standard complex structure on $ \mathbb{C}^n \cong \mathbb{R}^{2n}$ and the standard K\"ahler form on  $ \mathbb{C}^n \cong \mathbb{R}^{2n}$ respectively.  The closed complex-valued $n$-form is given by
\[
\Omega = dz_1\wedge \cdot\cdot\cdot \wedge dz_n,
\]
where $z_j = x_j+ i y_j$ are complex coordinates of $\mathbb{C}^n$.

An $n$-dimensional oriented submanifold $\Sigma$ immersed into $ \mathbb{C}^n $ is called a {\it Lagrangian submanifold} if $\left. \omega \right|_{\Sigma}=0$.  For a Lagrangian submanifold $\Sigma$,  the induced volume form $d\mu_{\Sigma}$ on $\Sigma$ and the complex-valued $n$-form $\Omega$ are  related by 
\[
\left. \Omega \right|_{\Sigma} = e^{i\theta}d\mu_{\Sigma} = \cos\theta d\mu_{\Sigma} + i \sin\theta d\mu_{\Sigma},
\]
where $\theta$ is some multi-valued function
called the {\it Lagrangian angle} and is well-defined up to an
additive constant $2k\pi, k\in \mathbb{Z}$. Nevertheless, $\cos\theta$ and $\sin\theta$ are single valued functions on 
$\Sigma$. If $\cos\theta > 0$, then $\Sigma$ is said to be {\it almost calibrated}. The relation between the Lagrangian angle and
the mean curvature is given by (see \cite{SW})
 \begin{equation}\label{eqn-HJL}\aligned
H=J\n \theta.
\endaligned
\end{equation}

It is known that one can study the existence of minimal Larangian submanifolds via the mean curvature flow (c.f. \cite{TY02}).

Let $X: M^{n} \rightarrow \mathbb{R}^{m+n}$ be an isometric immersion from an $n$-dimensional oriented Riemannian manifold $M$ to the Euclidean space $\mathbb{R}^{m+n}$.  
The mean curvature flow (MCF) in Euclidean space is a one-parameter family of immersions $X_t= X(\cdot, t): M^n \rightarrow \mathbb{R}^{m+n}$ with the corresponding images $M_t=X_t(M)$ such that
\begin{equation}\label{eqn-MCF}
\begin{cases}\aligned
\frac{\partial}{\partial t}X(x,t)=& H(x,t), x\in M,\\
X(x,0)=&X(x),
\endaligned
\end{cases}
\end{equation}
is satisfied, where $H(x, t)$ is the mean curvature vector of $M_t$ at $X(x, t)$ in $\mathbb{R}^{m+n}$.

Smoczyk \cite{Smo96} showed that the Lagrangian property is preserved along the mean curvature flow, in which case the flow is called the Lagrangian mean curvature flow. By establishing a new monotonicity formula, Chen--Li \cite{CL01} and Wang \cite{Wan01} independently proved that there is no Type ${\rm I}$ singularity for an almost calibrated Lagrangian mean curvature flow.  Later, Neves \cite{Nev07} demonstrated that the same conclusion holds for zero-Maslov class. Therefore the geometric and analytic nature of Type
II singularities in the Lagrangian mean curvature flow with zero-Maslov class has been of great interest.

One of the most important examples of Type II singularities is the translating soliton. Recall that $\Sigma^n$ is called a {\it translating soliton} in $\mathbb{R}^{2n}$ if it
 satisfies
\begin{equation}\label{eqn-T111}
H= V_{0}^N,
\end{equation}
where  $V_0$ is a fixed vector in $\mathbb{R}^{2n}$ with unit length and $V_{0}^N$ denotes the orthogonal projection of $V_0$ onto the normal bundle of $\Sigma^n$.
There are plenty of works on the subject of Lagrangian translating solitons, see \cite{CL12, CL15, CheQiu16, Coo15, LLQ21, LLS24,Qiu23, Smo19, Sun13, Sun14, Sun15,  SW} and the references therein. Joyce-Lee-Tsui \cite{JLT10} constructed lots of Lagrangian
translating solitons with oscillation of the Lagrangian
angle arbitrarily small, they are important in studying the regularity of Lagrangian MCF. In \cite{HS10}, Han-Sun showed that if the cosine of the Lagrangian angle has a positive lower bound, then any complete Lagrangian translating soliton with nonnegative sectional curvature in $\mathbb{R}^4$ has to be an affine plane.
Later, Neves-Tian \cite{NT} proved that under some natural conditions, if the first Betti number of the Lagrangian translating soliton $\Sigma^2 \subset \mathbb R^4$ is finite, and either $\int_{\Sigma^2} |H|^2 d\mu$ is finite or $\Sigma^2$ is static and almost calibrated, then $\Sigma^2$ is a plane. The restriction on the dimension of the ambient space $\mathbb{R}^4$ is necessary (see \cite{NT}). Inspired by the previous results, it is natural to study that whether we can extend the corresponding result of \cite{HS10} to higher dimension. Along this direction, the author and Zhu \cite{QZ22} obtained that if the cosine of the Lagrangian angle has a positive lower bound, then any proper Lagrangian translating soliton with nonnegative scalar curvature in $\mathbb{R}^{2n}$ has to be an affine $n$-plane.

In this paper, we continue to study the rigidity of Lagrangian translating solitons. With the aid of appropriate auxiliary functions, we obtain point-wise estimates for the mean curvature, as shown by (\ref{eqn-H3}). This leads to the following rigidity result of Lagrangian translating solitons:

\begin{thm}\label{thm-lag}

Let $X: \Sigma^n \to \mathbb{R}^{2n}$ be a complete $n$-dimensional Lagrangian translating soliton with nonnegative scalar curvature. Assume that $\cos\theta>0$ and
\begin{equation}\label{eqn-thm11}
\lim_{a\to +\infty} \frac{\sup_{B_a(x_0)}(\cos\theta)^{-1}}{\sqrt{a}} = 0,
\end{equation}
where $B_a(x_0)$ is a closed geodesic ball of $\Sigma^n$ center at $x_0$ of radius $a$ and $\theta$ is the Lagrangian angle.
Then $\Sigma^n$ has to be an affine $n$-plane.

\end{thm}

When $\cos\theta$ has a positive lower bound, then it satisfies the growth condition (\ref{eqn-thm11}) in Theorem \ref{thm-lag}. Hence we obtain 

\begin{cor}\label{cor-lag}

Let $X: \Sigma^n \to \mathbb{R}^{2n}$ be a complete $n$-dimensional Lagrangian translating soliton with nonnegative scalar curvature and $\cos\theta \geq \d$, where $\theta$ is the Lagrangian angle and $\d$ is a positive constant.
Then $\Sigma^n$ has to be an affine $n$-plane.

\end{cor}

\begin{rem}

(1). While the ambient space is the four dimensional Euclidean space $\mathbb R^4$, Han-Sun \cite{HS10} proved that any complete Lagrangian translating soliton in $\mathbb{R}^4$ 
with nonnegative sectional curvature and $\cos\theta \geq \d$ must be an affine plane. 

(2). By using a different method, the author and Zhu \cite{QZ22} showed that the above Corollary \ref{cor-lag} holds under an additional properness assumption. Hence the condition in Theorem \ref{thm-lag} is much weaker than the one in \cite{QZ22}. 
\end{rem}

\begin{rem}

The assumption $\cos\theta > 0$ in Theorem \ref{thm-lag} is necessary. Notice that the ``grim reaper" $(-\ln\cos x, x, y, 0), |x|< \frac{\pi}{2}, y\in \mathbb{R}$, is a Lagrangian translating soliton which translates in the direction of the constant vector $(1, 0, 0, 0)$. By direct computation, we have $\cos\theta = 0$ for the  ``grim reaper".

\end{rem}

Another important class of translating solitons in higher codimension is the symplectic translating solitons, which are special solutions to symplectic mean curvature flows (see \cite{CL01, Wan01}). Symplectic translating solitons play an important role in analysis of singularities in symplectic MCFs. Let $\o, \la \cdot, \cdot \ra$ be the standard K\"ahler form and Euclidean metric on $\mathbb{R}^4$, respectively. The {\it K\"ahler angle} $\a$ of $\Sigma^2$ in $\mathbb{R}^4$ is defined by 
\[
\left. \o \right|_{\Sigma}=\cos\a d\mu_{\Sigma},
\]  
where $d\mu_{\Sigma}$ is the area element of the induced metric. 

A surface $\Sigma^2$ is called {\it symplectic} if the K\"ahler angle satisfies $\cos\a>0$. By the formula (2.5) in \cite{HS10} (see also (4.4) in \cite{Qiu22b}),
\begin{equation}\label{eqn-ka8}\aligned
|\overline \n J_\Sigma|^2 = |B|^2 - 2R_{1234},
\endaligned
\end{equation}
where $R_{1234}$ is the normal curvature of $\Sigma^2$ and $|\overline \n J_\Sigma|^2= |h^2_{11}+h^1_{12}|^2 + |h^2_{21}+h^1_{22}|^2 + |h^2_{12}-h^1_{11}|^2 + |h^2_{22}-h^1_{21}|^2$ which depends only on the orientation of $\Sigma^2$  and does not depend on the choice of local coordinates.

 Han–Li \cite{HL09} gave an estimate of the kähler angle for symplectic standard translating solitons by using the Omori-Yau maximum principle. Afterward, Han-Sun \cite{HS10} proved that if the cosine of the Kähler
angle has a positive lower bound, then any complete symplectic translating soliton
with nonpositive normal curvature and bounded second fundamental form has to be
an affine plane. See also the related work by the author \cite{Qiu22b}.

By using a similar approach of the proof of Theorem \ref{thm-lag},  we derive a rigidity result of symplectic translating solitons as follows:

\begin{thm}\label{thm-sym}

Let $X: \Sigma^2 \to \mathbb{R}^4$ be a complete symplectic translating soliton with nonpositive normal curvature. Suppose that the second fundamental form is bounded by a positive constant $C_0$. Assume that
\[
\lim_{a\to +\infty} \frac{\sup_{B_a(x_0)}(\cos\a)^{-1}}{\sqrt{a}} = 0,
\]
where $B_a(x_0)$ is a closed geodesic ball of $\Sigma^n$ center at $x_0$ of radius $a$ and $\a$ is the K\"ahler angle. Then $\Sigma^2$ has to be an affine plane.
\end{thm}

\begin{rem}

Han-Sun \cite{HS10} showed that if $\cos\a$ has a positive lower bound, then any complete symplectic translating soliton with bounded second fundamental form and nonpositive normal curvature has to be an affine plane. It follows that $(\cos\a)^{-1}$ is bounded. Hence the condition in Theorem \ref{thm-sym} is weaker than the one in \cite{HS10}.

\end{rem}

In the proofs of the above Theorems \ref{thm-lag} and \ref{thm-sym} we take $r$ as the intrinsic distance of the translator, then by using the $V$-Laplacian comparison theorem (see Theorem 3 in \cite{CJQ12}), we can estimate $\mathcal{L}_{II} r:= \D r +\la V, \n r \ra$, where $V:= V_0^{T}$.  If we take $r:= |X|$ as the extrinsic distance instead of the intrinsic distance, we can also obtain 

\begin{thm}\label{thm-sym2}

Let $X: \Sigma^2 \to \mathbb{R}^4$ be a proper symplectic translating soliton with nonpositive normal curvature.  
Assume that
\[
(\cos\a)^{-1} = o(\sqrt{R}),
\]
where $R$ is the Euclidean distance from any point in $\Sigma^2$.
 Then $\Sigma^2$ has to be an affine plane.
\end{thm}

\begin{rem}

The author \cite{Qiu22b} proved that if $\cos \a$ has a positive lower bound, then any proper symplectic translating soliton with nonpositive normal curvature has to be an affine plane. In this case, $(\cos\a)^{-1}$ is bounded. Hence Theorem \ref{thm-sym2} improves the corresponding result of \cite{Qiu22b}. 

\end{rem}

\vskip24pt

\section{Notations}

\bigskip

Let $\Sigma^n$ be a submanifold in $\mathbb{R}^{2n}$. Let $\overline{\n}, \n$ be the Levi-Civita connection on $\mathbb R^{2n}$ and $\Sigma^n$, respectively.
The second fundamental form $B$ of $\Sigma^{n}$ in $\mathbb{R}^{2n}$ is defined by
\[
B_{UW}:= (\overline{\n}_U W)^N
\]
for $U, W \in \Gamma(T\Sigma^n)$. We use the notation $( \cdot )^T$ and $(
\cdot )^N$ for the orthogonal  projections into the tangent bundle
$T\Sigma^n$ and the normal bundle $N\Sigma^n$, respectively. For $\nu \in
\Gamma(N\Sigma^n)$ we define the shape operator $A^\nu: T\Sigma^n \rightarrow T\Sigma^n$
by
\[
A^\nu (U):= - (\overline{\n}_U \nu)^T.
\]
Taking the trace of $B$ gives the mean curvature vector $H$ of $\Sigma^n$
in $\mathbb{R}^{2n}$ and
\[
H:= \hbox{trace} (B) = \sum_{i=1}^{n} B_{e_ie_i},
\]
where $\{ e_i \}$ is a local orthonormal frame field of $\Sigma^n$.

Let $V$ be a smooth vector field on $\Sigma^n$ and $g$ the metric of $\Sigma^n$. The Bakry--Emery Ricci curvature is defined by 
\[
{\rm Ric}_V := {\rm Ric} - \frac{1}{2}L_V g,
\]
where ${\rm Ric}$ is the Ricci curvature of $\Sigma^n$ and $L_V$ is the Lie derivative.

\vskip24pt

\section{Proofs of main theorems}

\vskip8pt

Let  $\D$ be the Laplace-Beltrami operator on $\Sigma^n$ and $V:=V_0^{T}$. Denote $\mathcal{L}_{II}:= \D + \la V, \n\cdot \ra$.

\vskip12pt

\noindent{\bf Proof of Theorem \ref{thm-lag}}
Fix a point $p \in \Sigma^n$ and consider an orthonormal frame $\{e_i\}$, defined on a neighborhood around $p$ such that $\n e_i|_{p} = 0$. By (\ref{eqn-HJL}), we get
 \begin{equation*}\label{eqn-HJ}\aligned
e_i(\theta) = \langle H, Je_i \rangle.
\endaligned
\end{equation*}
It follows that
 \begin{equation}\label{eqn-LLA}\aligned
\D \theta = \sum_{i=1}^n e_i\langle H, Je_i \rangle = \sum_{i=1}^n \langle \overline{\n}_{e_i}H, Je_i \rangle + \sum_{i=1}^n \langle H, \overline{\n}_{e_i}(Je_i) \rangle.
\endaligned
\end{equation}
Since $\Sigma$ is Lagrangian, we derive
 \begin{equation}\label{eqn-LLA1}\aligned
 \sum_{i=1}^n \langle H, \overline{\n}_{e_i}(Je_i) \rangle = \langle H, JH \rangle = 0.
\endaligned
\end{equation}
The translating soliton equation (\ref{eqn-T111}) and (\ref{eqn-HJL}) imply that
 \begin{equation}\label{eqn-LLA2}\aligned
\sum_{i=1}^n \langle \overline{\n}_{e_i}H, Je_i \rangle = & \sum_{i=1}^n \langle \overline{\n}_{e_i}(V_0 - V), Je_i \rangle = - \sum_{i=1}^n \langle \overline{\n}_{e_i}V, Je_i \rangle \\
=& - \sum_{i=1}^n e_i \langle V, Je_i \rangle +  \sum_{i=1}^n\langle V, \overline{\n}_{e_i}(Je_i) \rangle \\
=&  \sum_{i=1}^n \langle V, J\overline{\n}_{e_i}e_i \rangle = \langle V, JH \rangle = - \langle V, \n \theta \rangle.
\endaligned
\end{equation}
Substituting (\ref{eqn-LLA1}) and (\ref{eqn-LLA2}) into (\ref{eqn-LLA}), we obtain 
 \begin{equation}\label{eqn-LA}\aligned
\D \theta = -\langle V, \n \theta \rangle.
\endaligned
\end{equation}
By (\ref{eqn-LA}) and (\ref{eqn-HJL}), we have
 \begin{equation*}\label{eqn-CLA}\aligned
\D \cos\theta =& -\cos\theta|\n \theta|^2-\sin\theta\D\theta \\
= &-\cos\theta|\n \theta|^2+\sin\theta\langle V, \n \theta \rangle \\
=& -\cos\theta|H|^2- \langle V, \n \cos \theta \rangle.
\endaligned
\end{equation*}
Namely
\begin{equation}\label{eqn-la1}\aligned
\mathcal{L}_{II} \cos \theta = -\cos\theta |H|^2.
\endaligned
\end{equation}
Let 
\begin{equation}\label{eqn-la2}\aligned
h_1= (\cos \theta)^{-2}.
\endaligned
\end{equation}
It follows that
\begin{equation}\label{eqn-la3}\aligned
\mathcal{L}_{II} h_1 = & - 2 (\cos \theta)^{-3}\mathcal{L}_{II} \cos\theta + 6 |\n(\cos \theta)^{-1} |^2 \\
=& 2(\cos \theta)^{-2} |H|^2 + \frac{3}{2}(\cos \theta)^{2} |\n (\cos \theta)^{-2}|^2 \\
=& 2h_1 |H|^2 + \frac{3}{2}h_1^{-1}|\n h_1|^2.
\endaligned
\end{equation}
 From the translating soliton equation (\ref{eqn-T111}), we derive
\begin{equation*}
\n_{e_j} H =  \left( \overline{\n}_{e_j}(V_0-\la V_0, e_k \ra e_k) \right)^N = - \la V_0, e_k \ra B_{e_j e_k}
\end{equation*}
and
\begin{equation*}
\n_{e_i}\n_{e_j} H = - \la V_0, e_k \ra \n_{e_i} B_{e_j e_k} -\la H, B_{e_i e_k} \ra B_{e_j e_k}.
\end{equation*}
Hence using the Codazzi equation, we obtain that
 \begin{equation*}\aligned
\mathcal{L}_{II} |H|^2 = & \D |H|^2 +\la V, \n |H|^2 \ra  \\
=& 2\la \n_{e_i}\n_{e_i}H, H \ra +2|\n H|^2+ \la V, \n |H|^2\ra \\
=& -2 \la H, B_{e_i e_k} \ra^2 - 2\la \n_{V_0^{T}} H, H \ra +2|\n H|^2 + \la V, \n |H|^2\ra \\
=& -2 \la H, B_{e_i e_k} \ra^2 - \n_V |H|^2 +2|\n H|^2 + \la V, \n |H|^2\ra \\
=&  -2 \la H, B_{e_i e_k} \ra^2  + 2|\n H|^2.
\endaligned
\end{equation*}
It follows that 
 \begin{equation}\label{eqn-la4}\aligned
\mathcal{L}_{II} |H|^2 \geq   2|\n H|^2-2|B|^2|H|^2\geq   2|\n |H||^2-2|B|^2|H|^2.
\endaligned
\end{equation}
By the Gauss equation, we get
\begin{equation}\label{eqn-Gauss}
\widetilde R(e_i, e_j, e_i, e_j) = \left< B_{e_i e_i}, B_{e_j e_j} \right> - \left< B_{e_i e_j}, B_{e_j e_i} \right>,
\end{equation}
where $\widetilde R$ is the curvature operator of $\Sigma$. Thus the scalar curvature $S$ of $\Sigma$ satisfies
 \begin{equation}\label{eqn-Sca}\aligned
S= \sum_{i,j} \widetilde R(e_i, e_j, e_i, e_j) =  |H|^2 - |B|^2.
\endaligned
\end{equation}
Combining (\ref{eqn-la4}) with (\ref{eqn-Sca}), it follows 
\begin{equation}\label{eqn-la5}\aligned
\mathcal{L}_{II} |H|^2 \geq 2|\n|H||^2 - 2|H|^4+2S|H|^2.
\endaligned
\end{equation}
 Let $\varphi = \varphi( h_{1})$ be a smooth nonnegative and nondecreasing function of $ h_{1}$ to be determined later. Then by (\ref{eqn-la3}) we have
 \begin{equation}\label{eqn-la6}\aligned
\mathcal{L}_{II} \varphi = & \D \varphi + \la V, \n \varphi \ra = \varphi'' |\n h_1|^2 + \varphi' \mathcal{L}_{II}  h_1 \\
= & \varphi'' |\n  h_1|^2 + \varphi'\left( 2h_1|H|^2+ \frac{3}{2}h_1^{-1}|\n h_1|^2 \right).
\endaligned
\end{equation}
From (\ref{eqn-la5}) and (\ref{eqn-la6}), we derive
 \begin{equation}\label{eqn-la7}\aligned
\mathcal{L}_{II} (|H|^2\varphi) = & (\mathcal{L}_{II} |H|^2 )\varphi + 2\la \n|H|^2, \n \varphi \ra + |H|^2 \mathcal{L}_{II}\varphi \\
\geq & \left( 2|\n |H||^2-2|H|^4 +2S|H|^2\right)\varphi +2 \la \n |H|^2, \n \varphi \ra\\
&+|H|^2 \left( \varphi'' |\n  h_1|^2 + \varphi'\left( 2h_1|H|^2+ \frac{3}{2}h_1^{-1}|\n h_1|^2 \right) \right).
\endaligned
\end{equation}
Since 
\begin{equation}\label{eqn-la55}\aligned
2\la \n|H|^2, \n \varphi \ra = &\varphi^{-1}\la \n(|H|^2\varphi), \n \varphi \ra - \varphi^{-1}|H|^2|\n \varphi|^2 + 2|H|\la \n|H|, \n \varphi \ra  \\
\geq &\varphi^{-1}\la \n(|H|^2\varphi), \n \varphi \ra - \varphi^{-1}|H|^2|\n \varphi|^2  \\
&-2|\n|H||^2\varphi -\frac{1}{2}|H|^2|\n\varphi|^2\varphi^{-1} \\
= &\varphi^{-1}\la \n(|H|^2\varphi), \n \varphi \ra -2|\n|H||^2\varphi -\frac{3}{2}|H|^2|\n\varphi|^2\varphi^{-1}.
\endaligned
\end{equation}
Therefore by (\ref{eqn-la7}), (\ref{eqn-la55}) and the scalar curvature $S\geq 0$, we have
 \begin{equation}\label{eqn-la77}\aligned
\mathcal{L}_{II} (|H|^2\varphi) \geq & \left( \varphi'' +\frac{3}{2}\varphi'  h_1^{-1} - \frac{3}{2}\varphi'^2\varphi^{-1} \right)|H|^2|\n h_1|^2 \\
& + \left( 2\varphi'h_1 - 2\varphi \right)|H|^4 + \varphi^{-1}\la \n(|H|^2\varphi), \n \varphi \ra.
\endaligned
\end{equation}
Let 
\[
k_1= \frac{1}{2}\inf_{B_a(x_0)} h_1^{-1}, \quad \varphi=\frac{ h_1}{1-k_1 h_1}.
\] 
We know that $k_1$ is a positive constant due to the fact that $B_a(o)$ is compact.
Then by direct computation, we have
\begin{equation*}\label{eqn-var1}\aligned
\varphi'= \frac{1}{(1-k_1 h_1)^2}, \quad \quad \varphi''=\frac{2k_1}{(1-k_1 h_1)^3},
\endaligned
\end{equation*}
\begin{equation}\label{eqn-var111}\aligned
&\varphi'' + \frac{3}{2}\varphi' h_1^{-1} - \frac{3}{2} \varphi'^2\varphi^{-1} \\
= & \frac{2k_1}{(1-k_1 h_1)^3} + \frac{3}{2} \left( \frac{1}{ h_1(1-k_1 h_1)^2}  -  \frac{1}{(1-k_1 h_1)^4} \cdot \frac{1-k_1 h_1}{h_1} \right) \\
=& \frac{2k_1}{(1-k_1 h_1)^3} + \frac{3}{2} \left( \frac{1}{h_1(1-k_1 h_1)^2}  -  \frac{1}{ h_1 (1-k_1 h_1)^3}  \right) \\
=& \frac{k_1}{2(1-k_1 h_1)^3} = \frac{k_1}{2 h_1(1-k_1 h_1)^2}\varphi,
\endaligned
\end{equation}
\begin{equation}\label{eqn-var222}\aligned
2\varphi' h_1 - 2\varphi = & \frac{2 h_1}{(1-k_1 h_1)^2} -\frac{2 h_1}{1-k_1 h_1} 
= \frac{2k_1 h_1^{2}}{(1-k_1 h_1)^2} = 2k_1\varphi^2.
\endaligned
\end{equation}
\begin{equation}\label{eqn-var333}\aligned
\varphi^{-1}\n \varphi = \frac{1-k_1 h_1}{ h_1} \cdot \frac{1}{(1-k_1 h_1)^2} \n  h_1=\frac{1}{ h_1 (1-k_1 h_1)}\n h_1.
\endaligned
\end{equation}
Substituting (\ref{eqn-var111})--(\ref{eqn-var333}) into (\ref{eqn-la77}), we obtain 
\begin{equation}\label{eqn-la12}\aligned
\mathcal{L}_{II}(|H|^2\varphi) \geq &   2k_1|H|^4\varphi^2 + \frac{k_1}{2h_1(1-k_1 h_1)^2} |\n  h_1|^2\cdot |H|^2 \varphi \\
&+ \frac{1}{ h_1(1-k_1 h_1)} \la \n  h_1, \n(|H|^2\varphi)\ra.
\endaligned
\end{equation}
From the equation (\ref{eqn-Gauss}), we know that
  \begin{equation}\label{eqn-Ric}\aligned
{\rm Ric}(e_i, e_i) = & \sum_{j}\left< B_{e_i e_i}, B_{e_j e_j} \right> - \sum_j\left< B_{e_i e_j}, B_{e_j e_i} \right> \\
=& H^\a h^\a_{ii} - \sum_j (h^\a_{ij})^2,
\endaligned
\end{equation}
where $h^\a_{ij}=\la B_{e_i e_j}, \nu_\a \ra, H^\a =\sum_j h^\a_{jj}$, and $\{\nu_\a\}$ is a local orthonormal frame field on $N\Sigma$.  Since 
 \begin{equation}\label{eqn-Lie1}\aligned
L_V e_i =& \overline\n_Ve_i - \overline\n_{e_i}V = \la V_0, e_j \ra \overline\n_{e_j}e_i - \overline\n_{e_i}(\la V_0, e_j \ra e_j) \\
=& \la V_0, e_j \ra B_{e_i e_j} - \la V_0, B_{e_i e_j} \ra e_j - \la V_0, e_j \ra B_{e_i e_j} \\
=& - \la V_0, B_{e_i e_j} \ra e_j = - \la H, B_{e_i e_j} \ra e_j = - H^\a h^\a_{ij}e_j.
\endaligned
\end{equation}
Therefore we get
 \begin{equation}\label{eqn-Lie2}\aligned
-\frac{1}{2}(L_V g)(e_i, e_i) =& -\frac{1}{2}\left( V\la e_i, e_i \ra -2g(L_V e_i, e_i) \right) \\
= &g(L_V e_i, e_i) = \la - H^\a h^\a_{ij}e_j, e_i \ra = -H^\a h^\a_{ii}.
\endaligned
\end{equation}
By (\ref{eqn-Ric}) and (\ref{eqn-Lie2}), we obtain
 \begin{equation}\label{eqn-RicV}\aligned
{\rm Ric}_V (e_i, e_i) = & {\rm Ric}(e_i, e_i) - \frac{1}{2}(L_V g)(e_i, e_i) \\
=& H^\a h^\a_{ii} - \sum_j (h^\a_{ij})^2 -H^\a h^\a_{ii} \\
=&  - \sum_j (h^\a_{ij})^2\geq -|B|^2.
\endaligned
\end{equation}
Combining (\ref{eqn-Sca}), (\ref{eqn-RicV}) and the assumption that the scalar curvature $S\geq 0$, we derive
\begin{equation}\label{eqn-RicV2}\aligned
{\rm Ric}_V (e_i, e_i) \geq -|B|^2 \geq -|H|^2 = -|V_0^{N}|^2 \geq - |V_0|^2 = -1.
\endaligned
\end{equation}
Let $r$ be the distance function on $\Sigma^n$ from $x_0$. Hence by the $V$-Laplacian comparison theorem (see Theorem 3 in \cite{CJQ12}), we get
\begin{equation}\label{eqn-Lap1}\aligned
\mathcal{L}_{II} r \leq \sqrt{n-1} \coth\sqrt{\frac{1}{n-1}}r +1 \leq \sqrt{n-1} + \frac{n-1}{r} +1.
\endaligned
\end{equation}
This implies that
\begin{equation}\label{eqn-Lap2}\aligned
\mathcal{L}_{II} r ^2 = & 2r\mathcal{L}_{II} r + 2|\n r|^2 \\
\leq & 2r\left(\sqrt{n-1} + \frac{n-1}{r} +1 \right)+2 \\
=& 2\left(\sqrt{n-1}+1\right) r + 2n.
\endaligned
\end{equation}
Let $\eta=((a^2-r^2)_+)^2$. On the support of $\eta$,
\begin{equation}\label{eqn-Lap3}\aligned
\mathcal{L}_{II} \eta = & 2|\n r^2|^2 +2(a^2-r^2)\mathcal{L}_{II} (-r^2) \\
= & 8r^2|\n r|^2 -2(a^2-r^2) \mathcal{L}_{II} r^2.
\endaligned
\end{equation}
Combining (\ref{eqn-Lap2}) with (\ref{eqn-Lap3}), we obtain
\begin{equation}\label{eqn-Lap4}\aligned
\mathcal{L}_{II} \eta \geq 8r^2|\n r|^2 -2(a^2-r^2)\left( 2\left(\sqrt{n-1}+1\right) r + 2n \right).
\endaligned
\end{equation}
Since 
\begin{equation}\label{eqn-grad-eta1}\aligned
2\la \n (|H|^2\varphi), \n \eta \ra = & 2\eta^{-1} \la \left( \n (|H|^2\varphi) \right)\eta, \n \eta \ra \\
=& 2\eta^{-1} \la \n(|H|^2\varphi\eta) -|H|^2\varphi\n \eta, \n\eta \ra \\
=& 2\eta^{-1} \la \n(|H|^2\varphi\eta), \n \eta \ra -32 |H|^2\varphi r^2 |\n r|^2.
\endaligned
\end{equation}
Therefore by (\ref{eqn-la12}), (\ref{eqn-Lap4}) and (\ref{eqn-grad-eta1}), we have
\begin{equation}\label{eqn-VLap-Bvareta1}\aligned
\mathcal{L}_{II} (|H|^2\varphi\eta) = & \left( \mathcal{L}_{II} (|H|^2\varphi) \right) \eta + 2\la \n(|H|^2\varphi), \n \eta \ra + |H|^2\varphi \mathcal{L}_{II}\eta\\
 \geq & 2k_1|H|^4\varphi^2\eta + \frac{k_1}{2 h_1(1-k_1 h_1)^2} |\n  h_1|^2 |H|^2 \varphi \eta  \\
& + \frac{1}{ h_1 (1-k_1 h_1)} \la \n  h_1, \n(|H|^2\varphi) \ra \eta \\
& + 2\eta^{-1} \la \n (|H|^2\varphi\eta), \n \eta \ra -24|H|^2\varphi r^2 |\n r|^2 \\
& - 2|H|^2\varphi (a^2-r^2)\left( 2\left(\sqrt{n-1}+1\right) r + 2n \right).
\endaligned
\end{equation}
The Cauchy inequality implies that
\begin{equation}\label{eqn-grad-eta2}\aligned
&\left| -\frac{1}{ h_1(1-k_1 h_1)}|H|^2\varphi \la \n h_1, \n \eta \ra \right| \\
\leq & \frac{1}{ h_1(1-k_1 h_1)}|H|^2\varphi |\n  h_1| |\n \eta| \\
= & \frac{1}{ h_1(1-k_1 h_1)}|H|^2\varphi |\n  h_1| \cdot 4\eta^{\frac{1}{2}}r|\n r| \\
\leq & \frac{k_1}{2 h_1(1-k_1 h_1)^2}|\n  h_1|^2 |H|^2 \varphi \eta + \frac{8}{k_1 h_1}|H|^2 \varphi r^2 |\n r|^2.
\endaligned
\end{equation}
It follows that
\begin{equation}\label{eqn-grad-eta3}\aligned
&\frac{1}{ h_1 (1-k_1 h_1)} \la \n  h_1, \n(|H|^2\varphi) \ra \eta \\
= &\frac{1}{ h_1 (1-k_1 h_1)} \left\la \n  h_1, \n(|H|^2\varphi \eta) \right \ra - \frac{1}{ h_1(1-k_1 h_1)}|H|^2\varphi \left\la \n h_1, \n \eta \right\ra \\
\geq & \frac{1}{ h_1 (1-k_1 h_1)} \la \n  h_1, \n(|H|^2\varphi \eta) \ra \\
& -  \frac{k_1}{2 h_1(1-k_1 h_1)^2}|\n  h_1|^2 |H|^2 \varphi \eta  - \frac{8}{k_1 h_1}|H|^2 \varphi r^2 |\n r|^2.
\endaligned
\end{equation}
Substituting (\ref{eqn-grad-eta3}) into (\ref{eqn-VLap-Bvareta1}), we obtain 
 \begin{equation}\label{eqn-VLap-Bvareta2}\aligned
\mathcal{L}_{II} (|H|^2\varphi\eta) \geq & 2k_1|H|^4\varphi^2\eta + \left\la 2\eta^{-1}\n \eta +\frac{\n h_1}{ h_1(1-k_1 h_1)}, \n(|H|^2\varphi \eta) \right\ra \\
& - 8\left( 3+\frac{1}{k_1 h_1} \right) |H|^2 \varphi r^2 |\n r|^2 \\
&- 2|H|^2 \varphi (a^2-r^2)\left( 2\left(\sqrt{n-1}+1\right) r + 2n \right).
\endaligned
\end{equation}
Let
\[
f= |H|^2\varphi \eta.
\]
Then the above inequality becomes
 \begin{equation}\label{eqn-VLap-Bvareta2}\aligned
\mathcal{L}_{II} f \geq & 2k_1\eta^{-1}f^2 + \left\la 2\eta^{-1}\n \eta +\frac{\n  h_1}{ h_1(1-k_1 h_1)}, \n f \right\ra \\
& - 8\left( 3+\frac{1}{k_1 h_1} \right)\eta^{-1} f r^2 |\n r|^2 - 2\eta^{-1}f (a^2-r^2)\left( 2\left(\sqrt{n-1}+1\right) r + 2n \right).
\endaligned
\end{equation}
Since $\left.f \right|_{\p_{B_a(x_0)}}=0$, $f$ achieves an absolute maximum in the interior of $B_a(x_0)$, say $f \leq f(q)$ for some $q$ inside $B_a(x_0)$.  We may assume $|H|(q)\neq 0$. By the maximum principle, we have
\[
\n f (q) = 0, \quad\quad \mathcal{L}_{II} f (q) \leq 0.
\]
Then by (\ref{eqn-VLap-Bvareta2}), we obtain the following at $q$:
 \begin{equation}\label{eqn-VLap-Bvareta3}\aligned
0\geq \mathcal{L}_{II} f \geq & 2k_1\eta^{-1}f^2 - 8\left( 3+\frac{1}{k_1 h_1} \right)\eta^{-1} f r^2 |\n r|^2 \\
& - 2\eta^{-1}f (a^2-r^2)\left( 2\left(\sqrt{n-1}+1\right) r + 2n \right).
\endaligned
\end{equation}
It follows that 
 \begin{equation}\label{eqn-B1}\aligned
f(q) \leq & \frac{1}{2k_1}\left( 8\left( 3+\frac{1}{k_1 h_1} \right) r^2 |\n r|^2 
+ 2(a^2-r^2)\left( 2\left(\sqrt{n-1}+1\right) r + 2n \right) \right) \\
\leq & \frac{1}{2k_1}\left( 8\left( 3+\frac{1}{k_1 } \right) a^2 
+ 2a^2\left( 2\left(\sqrt{n-1}+1\right) a + 2n \right) \right).
\endaligned
\end{equation}
For any $x \in \Sigma^n$, we can choose a sufficiently large $a$, such that $x\in B_{\frac{a}{2}}(x_0)$. 
 Thus from (\ref{eqn-B1}), we get
 \begin{equation}\label{eqn-H1}\aligned
|H|^2(x)\varphi(x)\eta(x) \leq  & \sup_{B_{\frac{a}{2}}(x_0)}f \leq  f(q) \\
\leq &  \frac{1}{2k_1}\left( 8\left( 3+\frac{1}{k_1 } \right) a^2 
+ 2a^2\left( 2\left(\sqrt{n-1}+1\right) a + 2n \right) \right).
\endaligned
\end{equation}
Since 
\[
\varphi=\frac{ h_1}{1-k_1 h_1}\geq \frac{1}{1-k_1}\geq 1, \quad\quad {\rm and} \quad \quad \eta\geq \frac{9}{16}a^4,
\]
Thus we obtain
\begin{equation}\label{eqn-H3}\aligned
|H|^2(x) \leq & \frac{8}{9a^4}\left( (24+4n)a^2+4\left(\sqrt{n-1}+1\right) a^3 \right)k_1^{-1} + \frac{8}{9a^4}\cdot 8a^2\cdot k_1^{-2} \\
\leq & C(n)\left(\left( \frac{1}{a} + \frac{1}{a^2} \right) \sup_{B_a(x_0)} h_1 + \frac{1}{a^2} \sup_{B_a(x_0)}h_{1}^2 \right).
\endaligned
\end{equation}
where $C(n)$ is a positive constant depending only on $n$.

Letting $a \to +\infty$ in (\ref{eqn-H3}), we then derive that $H\equiv 0$. Then by (\ref{eqn-Sca}) and the
assumption that the scalar curvature $S\geq 0$, we have $B\equiv 0$. Namely, $\Sigma^n$ is an affine $n$-plane.
\qed

\vskip12pt

\vskip12pt

\noindent{\bf Proof of Theorem \ref{thm-sym}}
By Proposition 3.1 in \cite{HL09}, we have
\begin{equation}\label{eqn-sym1}\aligned
\mathcal{L}_{II} \cos\a = -|\overline \n J_\Sigma|^2\cos\a.
\endaligned
\end{equation}
Direct computation gives us
 \begin{equation*}\label{eqn-ka3}\aligned
\mathcal{L}_{II} (\cos\a)^{-1} =& -  (\cos\a)^{-2}\mathcal{L}_{II}\cos\a + 2 (\cos\a)^{-3}|\n\cos\a |^2 \\
=& -  (\cos\a)^{-2}\mathcal{L}_{II}\cos\a + 2 \cos\a |\n  (\cos\a)^{-1}|^2.
\endaligned
\end{equation*}
It follows that
 \begin{equation}\label{eqn-ka4}\aligned
\mathcal{L}_{II}  (\cos\a)^{-2} =& 2 (\cos\a)^{-1}\mathcal{L}_{II} (\cos\a)^{-1} + 2|\n  (\cos\a)^{-1}|^2 \\
= & -2 (\cos\a)^{-3} \mathcal{L}_{II}\cos\a + 6|\n  (\cos\a)^{-1}|^2.
\endaligned
\end{equation}
From (\ref{eqn-sym1}) and (\ref{eqn-ka4}), we obtain
 \begin{equation}\label{eqn-ka5}\aligned
\mathcal{L}_{II}  (\cos\a)^{-2} =& 2 (\cos\a)^{-2} |\overline \n J_\Sigma|^2+ 6|\n  (\cos\a)^{-1}|^2 \\
=& 2 (\cos\a)^{-2} |\overline \n J_\Sigma|^2+ \frac{3}{2}(\cos\a)^{2}|\n (\cos\a)^{-2}|^2.
\endaligned
\end{equation}
Let
\[
h_2 = (\cos\a)^{-2}.
\]
Thus we get
\begin{equation}\label{eqn-ka6}\aligned
\mathcal{L}_{II} h_2 = 2h_2 |\overline \n J_\Sigma|^2+ \frac{3}{2}h_2^{-1}|\n h_2|^2.
\endaligned
\end{equation}
 Let $\varphi = \varphi( h_{2})$ be a smooth nonnegative and nondecreasing function of $ h_{2}$ to be determined later. Then we have
 \begin{equation}\label{eqn-ka77}\aligned
\mathcal{L}_{II} \varphi = & \D \varphi + \la V, \n \varphi \ra = \varphi'' |\n h_2|^2 + \varphi' \mathcal{L}_{II} h_2 \\
= & \varphi'' |\n  h_2|^2 + \varphi'\left( 2h_2 |\overline \n J_\Sigma|^2+ \frac{3}{2}h_2^{-1}|\n h_2|^2 \right).
\endaligned
\end{equation}
The inequalities (\ref{eqn-la4}) and (\ref{eqn-ka77}) imply that
 \begin{equation}\label{eqn-ka7}\aligned
\mathcal{L}_{II} (|H|^2\varphi) = & (\mathcal{L}_{II} |H|^2 )\varphi + 2\la \n|H|^2, \n \varphi \ra + |H|^2 \mathcal{L}_{II}\varphi \\
\geq & \left( 2|\n |H||^2-2|B|^2|H|^2 \right)\varphi +2 \la \n |H|^2, \n \varphi \ra\\
&+|H|^2 \left( \varphi'' |\n  h_2|^2 + \varphi'\left( 2h_2 |\overline \n J_\Sigma|^2+ \frac{3}{2}h_2^{-1}|\n h_2|^2 \right) \right).
\endaligned
\end{equation}
By (\ref{eqn-ka8}) and the assumption that the normal curvature $R_{1234} \leq 0$, we obtain
\begin{equation}\label{eqn-ka88}\aligned
|\overline \n J_\Sigma|^2 \geq |B|^2.
\endaligned
\end{equation}
From (\ref{eqn-la55}), (\ref{eqn-ka7}), (\ref{eqn-ka88}), we derive
\begin{equation}\label{eqn-ka9}\aligned
\mathcal{L}_{II} (|H|^2\varphi) 
\geq & \left( \varphi'' +\frac{3}{2}\varphi'  h_2^{-1} - \frac{3}{2}\varphi'^2\varphi^{-1} \right)|H|^2|\n h_2|^2 \\
& + \left( 2\varphi'h_2 - 2\varphi \right) |B|^2|H|^2 + \varphi^{-1}\la \n(|H|^2\varphi), \n \varphi \ra.
\endaligned
\end{equation}
By (\ref{eqn-RicV}) and the assumption that $|B| \leq C_0$, we have
\[
{\rm Ric}_V \geq -|B|^2g \geq -C_0^{2}g,
\]
where $g$ is the metric of $\Sigma^2$.

Let $r$ be the distance function on $\Sigma^2$ from $x_0$. The $V$-Laplacian comparison theorem implies that
\[
\mathcal{L}_{II} r \leq C_0 \coth\left(C_0r\right) +1 \leq  \frac{1}{r} + C_0 +1.
\]
It follows that
\begin{equation}\label{eqn-Lap22}\aligned
\mathcal{L}_{II} r ^2 = & 2r\mathcal{L}_{II} r + 2|\n r|^2 \\
\leq & 2r\left( \frac{1}{r} + C_0 +1 \right)+2 \\
=& 2\left(C_0+1\right) r + 4.
\endaligned
\end{equation}
 Let 
\[
k_2= \frac{1}{2}\inf_{B_a(x_0)} h_2^{-1}, \quad \varphi=\frac{ h_2}{1-k_2 h_2}, \quad \eta=((a^2-r^2)_+)^2.
\] 
Then by a similar proof of Theorem \ref{thm-lag}, we can conclude that for any $x\in \Sigma^2$, 
\begin{equation}\label{eqn-H33}\aligned
|H|^2(x) \leq & \frac{16}{9a^4}\left( 32a^2+4\left(C_0+1\right) a^3 \right)k_2^{-1} + \frac{16}{9a^4}\cdot 8a^2\cdot k_2^{-2} \\
\leq & C(C_0)\left(\left( \frac{1}{a} + \frac{1}{a^2} \right) \sup_{B_a(x_0)} h_2 + \frac{1}{a^2} \sup_{B_a(x_0)}h_{2}^2 \right),
\endaligned
\end{equation}
where $C(C_0)$ is a positive constant depending only on $C_0$.

Letting $a \to +\infty$ in the above inequality, we then obtain that $H\equiv 0.$ Then by Proposition 3.2 in \cite{HL09}, $B\equiv 0$. Namely, $\Sigma^2$ is an affine plane.
\qed

\vskip12pt

\noindent{\bf Proof of Theorem \ref{thm-sym2}}
Let $r=|X|$. Then we have
\begin{equation}\label{eqn-dist}\aligned
\n r^2 =& 2X^T, \quad \quad |\n r| \leq 1, \\
\mathcal{L}_{II} r^2 =& 4 + 2\la H, X \ra + \la V_0^{T}, 2X^T \ra \\
= & 4 + 2\la V_0^{N}, X \ra + 2\la V_0^{T}, X \ra \\
=& 4 + 2\la V_0, X\ra \leq 4+2r.
\endaligned
\end{equation}

Let $B_R(o)\subset \mathbb{R}^4$ be a closed ball of radius $R$ centered at the origin $o\in \mathbb{R}^4$ and $D_R(o):= \Sigma^2 \cap B_R(o)$.  Let
\[
\widetilde k_2= \frac{1}{2}\inf_{D_R(o)} h_2^{-1}, \quad \varphi=\frac{ h_2}{1-\widetilde k_2 h_2}, \quad \eta=((R^2-r^2)_+)^2.
\] 
Since $D_R(o)$ is compact by the fact that $\Sigma^2$ is proper, therefore $\widetilde k_2$ is a positive constant. Then by using the estimate (\ref{eqn-dist}) instead of (\ref{eqn-Lap22}) in the proof of Theorem \ref{thm-sym}, we can conclude that $\Sigma^2$ is an affine plane.
\qed

\begin{rem}

When the translating soliton $\Sigma^2$ is an entire graph in $\mathbb{R}^4$, 
 let $u: \Sigma \rightarrow \mathbb{C}$ be a smooth map and $\Sigma := \{ \left.(z, u(z))\right| z \in \mathbb{C} \}$, where $z$ is the complex coordinate of $\Sigma$. Let $u=u^1+\sqrt{-1}u^2$. By the work of Leung-Wan \cite{LeuWan07} (see also \cite{Qiu22b}), $\Sigma^2$ is hyper-Lagrangian with the complex phase map $\widetilde J: \Sigma^2 \to \mathbb{S}^2$.  In the following we will give the expression of this complex phase map:
 
  Choosing the canonical hyperk\"ahler structure on $\mathbb{C}^2$ such that the k\"ahler forms are
\begin{equation*}\label{ODE-add}
\begin{cases}\aligned
\o_1=& \frac{\sqrt{-1}}{2}\left( dz^1\wedge d\bar{z}^1 + dz^2\wedge d\bar{z}^2 \right) \\
\o_2 =& {\rm Re} \left( \sqrt{-1}dz^1\wedge dz^2 \right) \\
\o_3 =& {\rm Im}  \left( \sqrt{-1}dz^1\wedge dz^2 \right)
\endaligned
\end{cases}
\end{equation*}
and then $\O=\o_2 +\sqrt{-1}\o_3 = \sqrt{-1}dz^1\wedge dz^2 $ is the associated holomorphic symplectic 2-form on $\mathbb{C}^2$.
 It follows that
\[
\left. \o_1 \right|_{\Sigma}= \frac{\sqrt{-1}}{2}\left( 1+ |u_z|^2 - |u_{\bar z}|^2 \right) dz \wedge d{\bar z} = \left(1+ {\rm det}(du)\right)  \frac{\sqrt{-1}}{2}dz \wedge d{\bar z},
\]
where ${\rm det}(du)$ is the Jacobian of $u$, and
\[
\left. \O \right|_{\Sigma} = \sqrt{-1} u_{\bar z} dz\wedge d\bar z = 2u_{\bar z} \cdot  \frac{\sqrt{-1}}{2}dz\wedge d\bar z.
\]
 Let $\{\frac{\p}{\p x^1}, \frac{\p}{\p x^2}, \frac{\p}{\p x^3}, \frac{\p}{\p x^4}\}$ be the natural basis of $\mathbb{R}^4$. Let $e_1=\frac{\p}{\p x^1} + \frac{\p u^1}{\p x^1} \frac{\p}{\p x^3} + \frac{\p u^2}{\p x^1}\frac{\p}{\p x^4}, e_2=\frac{\p}{\p x^2} + \frac{\p u^1}{\p x^2} \frac{\p}{\p x^3} + \frac{\p u^2}{\p x^2}\frac{\p}{\p x^4}$. Hence the complex phase map can be represented by 
\[
\frac{1}{\sqrt{{\rm det}g}}\left( 1+ {\rm det}(du), {\rm Re}(2u_{\bar z}),  {\rm Im}(2u_{\bar z})  \right),
\]
where $g$ is the metric matrix, and $g_{ij}=\la e_i, e_j\ra$. Similarly, if we choose another hyperk\"ahler structure on $\mathbb{C}^2$ such that  the k\"ahler forms are
\begin{equation*}\label{ODE-add}
\begin{cases}\aligned
\o_1=& \frac{\sqrt{-1}}{2}\left( -dz^1\wedge d\bar{z}^1 + dz^2\wedge d\bar{z}^2 \right) \\
\o_2 =& {\rm Re} \left( -\sqrt{-1}d\bar{ z}^1\wedge dz^2 \right) \\
\o_3 =& {\rm Im}  \left( -\sqrt{-1}d\bar{z}^1\wedge dz^2 \right)
\endaligned
\end{cases}
\end{equation*}
and the complex phase map can be represented by 
\[
\frac{1}{\sqrt{{\rm det}g}}\left( {\rm det}(du) -1, {\rm Re}(2u_{z}),  {\rm Im}(2u_{z})  \right).
\]
Note that a surface in $\mathbb{R}^4$ being symplectic is equivalent to the condition that the image under the complex phase map is contained in an open hemisphere. Therefore by Theorem \ref{thm-sym2}, we have the following rigidity result:

Let $\Sigma^2 := \{ \left.(z, u(z))\right| z \in \mathbb{C} \}$ be an entire graphic translating soliton with nonpositive normal curvature, where $u: \Sigma \rightarrow \mathbb{C}$ is a smooth map.  Assume that either

(1) ${\rm det}(du) > -1$ and 
\[
\frac{\sqrt{{\rm det}g}}{1+{\rm det}(du)} = o(\sqrt{R})
\]
 or
 
 (2)  ${\rm det}(du) < 1$ and 
\[
\frac{\sqrt{{\rm det}g}}{1- {\rm det}(du)} = o(\sqrt{R}),
\]
where $R$ is the Euclidean distance from any point in $\Sigma^2$.
Then $\Sigma^2$ has to be an affine plane.

\end{rem}

 \vskip12pt

  \noindent{\bf Acknowledgements}    This work is partially supported by NSFC (No. 12471050) and Hubei Provincial Natural Science Foundation of China (No. 2024AFB746).

\vskip24pt

\end{document}